\documentclass[12pt,twoside,reqno]{amsart}

\usepackage{amsmath,amssymb,amsthm, amsfonts }
\usepackage{fancyhdr}
\usepackage{epsfig,color,colordvi}

\usepackage[numbers,sort&compress,square,comma]{natbib}

\theoremstyle{definition}

\newtheorem{example}{\sc Example}

\theoremstyle{remark}

\theoremstyle{plain}

\newcommand{\be}{\begin{equation}}
\newcommand{\ee}{\end{equation}}
\newcommand{\nn}{\nonumber}

\newcommand{\bb}{\mathbb{B}}

\def\ee{\mathbb{E}}
\def\nn{\mathbb{N}}
\def\rr{\mathbb{R}}
\def\pp{\mathbb{P}}

\def\BB{\mathcal B}

\def\XX{\mathcal X}

\newtheorem{thm}{Theorem}[section]
\newtheorem{cor}[thm]{Corollary}
\newtheorem{lem}[thm]{Lemma}
\newtheorem{prop}[thm]{Proposition}
\newtheorem{remarks}[thm]{Remarks}

\def\bdes{\begin{description}}
\def\ndes{\end{description}}
\def\beq{\begin{equation}}
\def\deq{\end{equation}}

\def\bdef{\begin{defn}}
\def\ndef{\end{defn}}

\def\bthm{\begin{thm}}
\def\nthm{\end{thm}}

\def\bprop{\begin{prop}}
\def\nprop{\end{prop}}

\def\brmk{\begin{remarks}}
\def\nrmk{\end{remarks}}

\def\bexa{\begin{exa}}
\def\nexa{\end{exa}}

\def\blem{\begin{lem}}
\def\nlem{\end{lem}}

\def\bcor{\begin{cor}}
\def\ncor{\end{cor}}

\def\dsp{\displaystyle}

\def\bexe{\begin{exe}}
\def\nexe{\end{exe}}

\def\bprf{\begin{proof}}
\def\nprf{\end{proof}}

\def\fQ{{{\rm Q}\kern-.65em {}^{{}_/ }\,}}
\def\fQQ{ {{\rm Q}\kern-.57em \scriptscriptstyle{}^{]\kern.055em[}\,}}

\def\ord{\kern0.1em o\kern-0.02em{}_{\ds\breve{}}\kern0.1em}
\def\Ord{\kern0.1em O\kern-0.02em{\ds\breve{}}\kern0.1em}
\def\ds{\displaystyle}

\def\fmonth{\ifcase\month\or Jan\or Feb\or Mar\or Apr
\or May\or Jun\or Jul\or Aug\or Sep
\or Oct\or Nov\or Dec\fi\ }
\def\mmddyyyy{\the\month.\the\day.\the\year}
\def\ddmmyyyy{\the\day.\the\month.\the\year}
\def\Mddyyyy{\fmonth~\the\day,~\the\year}

\providecommand{\pp}[1]{\langle#1\rangle}

\numberwithin{equation}{section}

\allowdisplaybreaks[1]

\begin{document}

%%%%%%%
%%% Author(s) and other info
%%%%%%%

\author{Yu Miao}   % First
\address{Department of Mathematics and Statistics, Wuhan
University, 430072 Hubei, China and College of Mathematics and
Information Science, Henan Normal University, 453007 Henan, China.}
%\curraddr{}
\email{yumiao728@yahoo.com.cn}
%\urladdr{www.math.ohio-state.edu/$\sim$firas}
%\thanks{NSF}
\author{Guangyu Yang}  % Second
\address{Department of Mathematics and Statistics,
Wuhan University, 430072 Hubei, China.}
%\curraddr{}
\email{study\_yang@yahoo.com.cn}
%\urladdr{www.math.wisc.edu/$\sim$seppalai}
%\thanks{NIH}
%\thanks{T.~Sepp\"al\"ainen was partially supported by National Science Foundation grant DMS-0402231.}

\date{December 26, 2005}
% this will put the date as a footnote on the first page

%\translator{Nobody}
%\dedicatory{Joe}


%%%%%%
%%% If you want an abstract
%%%%%%

%% Uncomment this if using \bfsection
%\renewcommand\abstractname{\noindent\bf Abstract}

%\frenchspacing

\begin{abstract}
 In the paper, the law of the iterated logarithm for additive functionals
of Markov
 chains is obtained under some weak conditions, which are weaker than
the conditions of
 invariance principle of additive functionals of Markov chains in M.
Maxwell and M. Woodroofe \cite{MMMW}
 (2000). The main technique is the martingale argument and the
 theory of fractional
coboundaries.

\end{abstract}

%%%%%%%
%%% Don't forget the title and date. To keep track of the notes.
%%% Preferably, name the file with something indicating the
%%% title and the date, as well.
%%%%%%%

%% You can copy the part bellow into the textfile

%% STARTING HERE

\title[The law of the iterated logarithm for Markov chains]
{The law of the iterated logarithm for additive functionals of
Markov chains}
% You can put the date manually or use \today, \mmddyyyy,
% \ddmmyyy, or \Mddyyyy  to have it automatically

\maketitle

%% and then define each to be what you wish (PP, fP, bbP, etc)
%\def\R{\fR}
%\def\N{\fN}
%\def\Z{\fZ}
%\def\Q{\fQ}
%\def\P{\fP}
%\def\E{\fE}
%\def\T{\fT}
%\def\C{\fC}

\def\X{\mathcal X}
\def\B{\mathcal B}

\section{Introduction}
Let $X_0, X_1, X_2,\ldots$ denote an ergodic stationary Markov chain
with values in a measurable space ($\XX,\BB$), transition function
$Q$, and stationary initial distribution $\pi$. Further, let
$L^2(\pi)$ denote the space of (equivalence classes of) square
integrable functions $g: \XX\to\rr$, for which $\|g\|^2:=\int_\XX
g^2d\pi<\infty$, and let $L^2_0(\pi)$ denote the set of $g\in
L^2(\pi)$ for which $\int_\XX gd\pi=0$. Given $g\in L^2_0(\pi)$, for
$n\geq 1$, let
 \beq\label{101}
 S_n=S_n(g):=\sum_{i=1}^n g(X_i).
 \deq
 To describe more detailedly the model of additive functionals of above
Markov chain,
 let $Q$ denote both the conditional distribution of
 $X_1$ given $X_0$ and the operator defined by
 $$Qh(x)=\int_\XX h(y)Q(x; dy)$$
 for a.e. $x\in\XX$ and all $h\in L^2(\pi)$. It is easy to see that
 $Q$ is a contraction. Let
 $$V_nh=\sum_{k=0}^{n-1}Q^kh$$
 for $h\in L^2(\pi)$ and $n\geq 1$.
Then each $V_n$ is a bounded linear operator. It is obviously that
$V_ng(x)=\ee[S_n(g)|X_1=x]$ for a.e. $x$ ($\pi$). For model
(\ref{101}), M. Maxwell and M. Woodroofe \cite{MMMW} (2000) gave the
central limit theorems under the following conditions:
 \beq\label{102}
g\in L^2_0(\pi),\ \  \sum_{n=1}^\infty n^{-3/2}\|V_n g\|<\infty,
 \deq
 where $\|\cdot\|$ denotes the norm in $L^2(\pi)$.

 It is well known that the law of the iterated
logarithm (in short LIL) is closely related to the central limit
theorem (in short CLT) in some sense. There are several approaches
to these
 problems. If the chain is Harris recurrent, then the problems may be
 reduced to the independent case in a certain sense, see S.P. Meyn and R.L. Tweedie
\cite{MT}
 (1993). If there is a solution to Poisson's equation, $h=g+Qh$,
 then the LIL and CLT problem may be reduced to the martingale case, see also
 M.I. Gordin and B.A. Lifsic \cite{GL} (1978) and S.P. Meyn and R.L.
Tweedie \cite{MT}
 (1993). R.N. Bhattacharya \cite{B} (1982) obtained CLT and LIL for
 ergodic stationary Markov process by discussion for infinitesimal
 generator.
 L.M. Wu \cite{Wu} (1999) extended the forward-backward martingale decomposition
of Meyer-Zheng-Lyons's type from the symmetric case to the general
stationary situation and gave Strassen's strong invariance
principle.

Our goal, in this paper, is to consider the problem that $S_n$
satisfies the LIL under some proper conditions. We will explore some
extensions of the Poisson's equation approach, along the lines of C.
Kipnis and S.R.S. Varadhan \cite{KV} (1986), and M. Maxwell and M.
Woodroofe \cite{MMMW} (2000), to the cases where a solution is not
required. We obtain the LIL for additive functionals of Markov
 chains in Section 2, where mainly depending on the LIL of martingale and the theory of
 fractional coboundaries developed by Y. Derriennic and M. Lin \cite{DL} (2001). In
Section 3, we will give two
 examples, whose central limit theorems and invariance principles
  are discussed in M. Maxwell and M. Woodroofe \cite{MMMW} (2000), to
illustrate our
  results.

\section{Main results}
For $\varepsilon>0$, let $h_\varepsilon$ be the solution to the
equation $(1+\varepsilon)h=Qh+g$,
 \beq\label{201}
 h_\varepsilon=\sum_{n=1}^\infty
\frac{Q^{n-1}g}{(1+\varepsilon)^n}=\varepsilon\sum_{n=1}^\infty
 \frac{V_{n}g}{(1+\varepsilon)^{n+1}}.
 \deq
 Let $\pi_1$ be the joint distribution of $X_0$ and $X_1$, so that
 $\pi_1(dx_0; dx_1)=Q(x_0; dx_1)\times \pi(dx_0)$; denote the
 norm in $L^2(\pi_1)$ by $\|\cdot\|_1$; and let
 $$H_\varepsilon(x_0, x_1)=h_\varepsilon(x_1)-Qh_\varepsilon(x_0)$$
 for $x_0$, $x_1\in\XX$. Then for any $\varepsilon>0$,
 $H_\varepsilon$ is in $L^2(\pi_1)$; the norm of $H_\varepsilon$ is
 $\|H_\varepsilon\|^2_1=\|h_\varepsilon\|^2-\|Qh_\varepsilon\|^2$;
 and $\int_\XX H_\varepsilon(x_0, x_1)Q(x_0; dx_1)=0$ for a.e. $x_0$
 $(\pi)$.

 Now, let us give a few more definitions. For $\varepsilon>0$,
 let
 $$M_n(\varepsilon)=H_\varepsilon(X_0,
X_1)+\cdots+H_\varepsilon(X_{n-1}, X_n)$$
 and
 $$R_n(\varepsilon)=Qh_\varepsilon(X_0)-Qh_\varepsilon(X_n),$$
 hence, by simple calculation,
 \beq\label{202}
 S_n(g)=M_n(\varepsilon)+\varepsilon
 S_n(h_\varepsilon)+R_n(\varepsilon).
 \deq
For any fixed $\varepsilon$, $M_n(\varepsilon)$ is a square
integrable martingale, where the martingale difference sequences is
stationary ergodic.

Before our main results, we need mention the following lemmas.

\begin{lem}\label{lem1}( See M. Maxwell and M. Woodroofe \cite{MMMW}
(2000) )
 \bdes
 \item{}{\rm ($ R1$)} If (\ref{102}) holds and $\|V_n
 g\|=O(n^\alpha)$ for some $\alpha>0$, then
 $\|h_\varepsilon\|=O(\varepsilon^{-\alpha})$ as $\varepsilon\to 0$.

 \item{}{\rm ($ R2$)} For $0<\varepsilon, \delta<\infty$,
 $\|H_\varepsilon-H_\delta\|_1^2\leq
 (\varepsilon+\delta)[ \|h_\varepsilon\|^2+\|h_\delta\|^2]$.

 \item{}{\rm ($ R3$)} If (\ref{102}) holds, then
$H=\lim_{\varepsilon\downarrow 0} H_\varepsilon$
 exists in $L^2(\pi_1)$.

 \item{}{\rm ($ R4$)} If (\ref{102}) holds, then
 \beq\label{203}
 S_n=M_n+R_n,
 \deq
where $M_1, M_2,\ldots$ and $R_1, R_2, \ldots$ have strictly
stationary increments, $M_1, M_2,\ldots$, is a square integrable
martingale, and $\ee (R_n^2)=o(n)$ as $n\to\infty$.
 \ndes
\end{lem}

Recall that $T$ is a Dunford-Schwartz (DS) operator on $L^1$ of a
probability space: if $T$ is a contraction of $L^1$ such that
$\|Tf\|_\infty\leq \|f\|_\infty$ for every $f\in L^\infty$. If
$\theta$ is a measure preserving transformation in a probability
space $(\Omega, \Sigma, \mu)$, then the operator $Tf=f\circ\theta$
is a DS operator on $L^1(\mu)$. More generally, any Markov operator
$\pp$ with an invariant probability yields a positive DS operator.

\begin{lem}\label{lem2} ( See Y. Derriennic and M. Lin \cite{DL}
(2001) )
 \bdes
 \item{{\rm $(L1)$}} Let $T$ be a contraction in a Banach space
 $X$, and let $0<\beta<1$. If $\dsp
 \sup_{n}\|\frac{1}{n^{1-\beta}}\sum_{k=1}^nT^ky\|<\infty$,
 then $y\in (I-T)^\alpha X$ for every $0<\alpha<\beta$.

 \item{{\rm $(L2)$}} Let $T$ is a DS operator in $L^1(\mu)$ of a
 probability space, and fix $1<p<\infty$, with dual $q=p/(p-1)$. Let
 $0<\alpha<1$, and $f\in(I-T)^\alpha L^p$. If
 $\alpha>1-\frac{1}{p}=\frac{1}{q}$, then
 $\dsp \frac{1}{n^{1/p}}\sum_{k=0}^{n-1}T^k f\to 0 $ a.e.
 \ndes
\end{lem}

\begin{lem}\label{lem3} ( See W.F. Stout \cite{WFS} (1970) )
Let ($Y_i, i\geq 1$) be a stationary ergodic stochastic sequence
with $\ee [Y_i|Y_1, Y_2, \cdots, Y_{i-1}]=0$ a.e. for all $i\geq 2$
and $\ee Y_1^2=1$. Then
$$\limsup_{n\to\infty}\frac{\sum_{i=1}^n Y_i}{(2n\log\log n)^{1/2}}=1,\
\  a.e.$$
\end{lem}

\begin{thm}\label{thm1} Let $g\in L^2_0(\pi)$ and
$\|V_ng\|=O(n^\alpha)$ for some $\alpha<1/2$. Then, we have
 \beq\label{204}
 \limsup_{n\to\infty}\frac{S_n}{(2n\log\log n)^{1/2}}=\|H\|_1,\ \  a.e.
 \deq
 where $H$ is defined as Lemma \ref{lem1} ($R3$).
\end{thm}
\begin{proof} At first, we notice that, by Lemma \ref{lem1} ($R4$),
$M_n$ is a square integrable martingale with strictly stationary
increments and, by Lemma \ref{lem3}, $M_n$ satisfies the law of the
iterated logarithm, i.e.
 \beq\label{205}
\limsup_{n\to\infty}\frac{M_n}{(2n\log\log n)^{1/2}}=\|H\|_1,\ \
a.e.
 \deq
where $\|H\|_1$ is due to the ergodic theorem, i.e.
$$\lim_{n\to\infty}\frac{1}{n}\sum_{i=1}^n H(X_{i-1}, X_i)^2=\|H\|_1^2,
\ \ a.e.$$

 Next we will consider $R_n$. For each $n\geq 1$, let
$\delta_j=2^{-j}$, $k_n$ be the unique integer $k$ for which
$2^{k-1}\leq n< 2^{k}$ and let $\varepsilon_n=2^{-k_n}$. Then, by
(\ref{202}) and Lemma \ref{lem1} ($R4$), we have
$R_n=M_n(\varepsilon_n)-M_n+\varepsilon_n
 S_n(h_{\varepsilon_n})+R_n(\varepsilon_n)$ and, therefore,
 \begin{eqnarray}\label{206}
 \ee(R^2_n)&\leq& 4\ee[M_n(\varepsilon_n)-M_n]^2+4\varepsilon^2_n\ee
 [S_n(h_{\varepsilon_n})]^2+4\ee [R_n(\varepsilon_n)]^2\nonumber\\
 &\leq&4n\|H_{\varepsilon_n}-H\|^2_1+4\|h_{\varepsilon_n}\|^2+8\|h_{\varepsilon_n}\|^2.
 \end{eqnarray}
By Lemma \ref{lem1} ($R1$), $\|h_{\varepsilon_n}\|=O(n^\alpha)$
 and by the definition of $\varepsilon_n$ and Lemma \ref{lem1} ($R2$),
there is a constant $C$
 for which
 \beq\label{207}
 \|H_{\varepsilon_n}-H\|_1\leq \sum_{j=k_n+1}^\infty
 \|H_{\delta_j}-H_{\delta_{j-1}}\|_1\leq C\sum_{j=k_n+1}^\infty
 \delta_j^{1/2-\alpha}=O(n^{\alpha-1/2}).
 \deq
 Thus, we have $\ee(R^2_n)=O(n^{2\alpha})$.

 On $\XX\times\XX$, define $f(x_0, x_1)=g(x_0)-H(x_0, x_1)$, then
 \beq\label{208}
 R_n=S_n-M_n=\sum_{i=1}^{n-1} [g(X_{i})-H(X_{i}, X_{i+1})]
=\sum_{i=1}^n f(X_{i}, X_{i+1}).
 \deq

  For sequence $x=(x_{i})_{i\in\nn}\in\XX^\nn$, define $F(x)=f(x_{0},
  x_{1})$ and let $\theta$ is the shift map on $\XX^\nn$ which is a
contraction on the space
  $L^2(\pi_1)$. Then, we have
\beq\label{209} F\in L^2(\pi_1) \;\;and\;\;
  R_n=\sum_{k=0}^{n-1}F\circ\theta^k.
  \deq
 From the fact $\ee(R^2_n)=O(n^{2\alpha})$, there exists a constant
$1/2<\beta<1-\alpha$, such that
 \beq\label{210}
\sup_{n}
\Big\|\frac{1}{n^{1-\beta}}\sum_{k=0}^{n-1}F\circ\theta^k\Big\|<\infty.
 \deq
Since Lemma \ref{lem2}
 ($L1$) and $\alpha<1/2$, we have $F\in(I-\theta)^\eta L^2(\pi_1)$,
 where $\eta\in(1/2, 1-\alpha)$. By Lemma \ref{lem2} ($L2$), we
 have
 $$\frac{1}{n^{1/2}}R_n\to 0,\ \ \pi_1-a.e.$$
Thus the result of the theorem holds by the fact, $S_n=M_n+R_n$.
\end{proof}
\begin{cor} If $g\in L^2(\pi)$, and either
 \beq\label{211}
 \sum_{k=1}^\infty k^{\delta-1/2}\|Q^k g\|<\infty
 \deq
 or
 \beq\label{212}
 \sum_{k=1}^\infty k^{\delta}\|Q^k g\|^2<\infty,
 \deq
 for some $\delta>0$, then (\ref{204}) holds.
\end{cor}
\begin{proof}
 By Kronecker's lemma and $V_ng=\sum_{k=0}^{n-1}Q^kg$, (\ref{211})
means that
 $$n^{\delta-1/2}\|V_ng\|\leq n^{\delta-1/2}\sum_{k=1}^n \|Q^k g\|\to
0, \ \ as\ \  n\to\infty.$$ It is obvious that, (\ref{212}) implies
 (\ref{211}). Form Theorem \ref{thm1}, the corollary holds.
\end{proof}

\begin{remarks}
{\rm If $g\in L^p_0(\pi)$, for some $p>2$, then (\ref{211}) or
(\ref{212}) are satisfied, M. Maxwell and M. Woodroofe \cite{MMMW}
(2000) gave an invariance principle of
$$\bb_n(t)=\frac{1}{\sqrt{n}}S_{\lceil nt\rceil},\ \ 0\leq t<1,$$
and $\bb_n(1)=\bb_n(1-)$ for $n\geq 1$, where $\lceil x\rceil$
denotes the least integer that is greater than $x$.}

\end{remarks}

\section{Some applications}
In this section, we will give two examples to show that they satisfy
the law of the iterated logarithm under some weakly dependent
conditions. Here, the processes do not have recurrent points and are
not strongly mixing. In fact, the two examples have been studied in
M. Maxwell and M. Woodroofe \cite{MMMW} (2000) and we should thank
their ideas and results.

\begin{example}{\rm \bf( Bernoulli shifts )}{\rm
Let $\varepsilon_k$, $k=0,\pm 1,\pm 2, \ldots$ be i.i.d. random
variables that take the values $0$ and $1$ with probability $1/2$
each and let
$$X_n=\sum_{k=1}^\infty(1/2)^{k+1}\varepsilon_{n-k}$$
for $n=0,1,2,\ldots$. Then $\{X_n$, $n\geq 0\}$, is an ergodic
stationary Markov chain taking values in $I=[0,1]$. The transition
function is defined by $Q(x;\{x/2\})=1/2=Q(x;\{1+x\}/2)$ for
$x\in[0,1]$ and the stationary initial distribution is the
restriction, $\lambda$ say, of Lebesgue measure to $I$. In this
case, $L^2(\lambda)$ is a familiar space, and it is easy to relate
the conditions in Theorem \ref{thm1} to regularity properties of
$g$.

 Proposition 3 in M. Maxwell and M. Woodroofe \cite{MMMW} (2000)
showed that if $g\in L^2_0(\lambda)$ and
 \beq\label{301}
\int_0^1\int_0^1[g(x)-g(y)]^2\frac{1}{|x-y|}\log^\delta\big[\frac{1}{|x-y|}\big]dxdy<\infty
 \deq
 for some $\delta>0$, then (\ref{212}) holds. For example, if
 $$g(x)=\frac{1}{x^\alpha}\sin\Big(\frac{1}{x}\Big),\ \ 0<x\leq 1,$$
 where $0<\alpha<1/2,$ then (\ref{212}) holds.}
 \end{example}
\begin{example}{\rm \bf( Lebesgue shifts )}{\rm
Let $U_k$, $k=0, \pm 1,\pm 2, \ldots$ be i.i.d. random variables
that are uniformly distributed over $I=[0,1]$ and let
$$X_n=(\cdots U_{n-2}, U_{n-1}, U_n),\ \ n\geq 0.$$
Then $\{X_n$, $n\geq 0\}$, is a stationary Markov process taking
values in $\XX=I^M$, where $M$ denote the non-positive integers. The
stationary initial distribution $\pi$ here is the countable product
of copies of Lebesgue measure. Processes of the form $g(X_n)$,
$n=0,1,2,\ldots$, include a wide class of stationary sequences.

 Define measures $\Gamma^1_\delta$ on $\XX\times\XX$ by
 $$\Gamma^1_\delta\{B\}=\sum_{k=1}^\infty
k^\delta(\pi\times\pi\times\lambda^k)\Big\{
 (x,y,z): \big[(x,z), (y,z)\big]\in B\Big\}$$
 for Borel set $B\subset\XX\times\XX$, where $x$ and ($x,z$) denote
 the sequences $x=(\ldots, u_{-2}, u_{-1}, u_0)$ and ($\ldots, u_{-1},
u_0, z_1, \ldots,
 z_k$) and $\lambda^k$ denotes Lebesgue measure on $I^k$.

 M. Maxwell and M. Woodroofe \cite{MMMW} (2000) showed that if
 $g\in L^2_0(\pi)$ and
 \beq\label{302}
 \int_\XX\int_\XX[g(x)-g(y)]^2\Gamma^1_\delta\{dx,dy\}<\infty,
 \deq
 for some $\delta>0$, then (\ref{212}) holds. To illustrate the
 (\ref{302}), observe that any $g\in L^2_{0}(\pi)$ may be written in
 the form
 $$g(x)=\sum_{k=0}^\infty g_k(u_{-k},\ldots, u_0),$$
 where $x=(\ldots, u_{-2}, u_{-1}, u_0)$, $g_k:\rr^{k+1}\to\rr$ are
 measurable, $\int_\rr g_k(u_{-k},\ldots,u_0) du_{-k}=0$ for
 a.e. $(u_{-k+1},\ldots,u_0)$, $\sum_{k=1}^\infty\int g^2_k
 d\lambda^{k+1}<\infty$, and $\lambda^k$ denotes $k-$dimensional
 Lebesgue measure.
 Then (\ref{212}) holds if for some $\delta>0$,
 $$\sum_{k=1}^\infty k^{1+\delta}\int g^2_k
 d\lambda^{k+1}<\infty.$$}
\end{example}

\section*{Acknowledgements}
The authors wish to thank Prof. L.M. Wu of Universit\'e Blaise
Pascal and Prof Y.P. Zhang of Wuhan University for their helpful
discussions and suggestions during writing this paper. And the
authors thank Prof. Y. Derriennic and Prof. M. Lin for sending the
key reference paper \cite{DL} to the authors.


\begin{thebibliography}{}

\bibitem{B} R.N. Bhattacharya, {\it On the functional central limit theorem and the law of the iterated logarithm
for Markov processes. } Z. Wahrs. verw. Gebiete., {\bf 60} (1982)
185-201.

\bibitem{CT} Y.S. Chow and H. Teicher, {\it Probability Theory.}
Sectond Edition. Springer-Verlag  New York 1988.

\bibitem{DL} Y. Derriennic and M. Lin, {\it Fractional Poisson
equations and ergodic theorems for fractional coboundaries.} Israel
J. math. {\bf 123} (2001), 93-130.

\bibitem{FE} W. Feller, {\it An introduction to probability theory and
its applications. Second Edition.} John Wiley and Sons, Inc. Vol 2,
1971.

\bibitem{GL} M.I. Gordin and B.A. Lifsic, {\it The central limit theorem for stationary Markov processes.}
Dokl. Akad. Nauk SSSR, {\bf 19} (1978), 392-394.


\bibitem{KV} C. Kipnis and S.R.S. Varadhan, {\it Central limit theorem for additive functionals of
reversible Markov processes and applications to simple exclusions.}
Comm. Math. Phys.,{\bf 104} (1986), 1-19.

\bibitem{MMMW} M. Maxwell and M. Woodroofe, {\it Central limit theorems
for additive functionals of Markov chains.} Ann. Probab. {\bf 28(2)}
(2000), 713-724.

\bibitem{MT} S.P. Meyn and R.L. Tweedie, {\it Markov chains and
stochastic stability.} Springer-Verlag  New York 1993.

\bibitem{PS} R.P. Pakshirajan and M. Sreehari, {\it The law of the
iterated logarithm for a Markov process.} Ann. Math. Statist. {\bf
41(3)} (1970), 945-955.

\bibitem{WFS} W.F. Stout, {\it The Hartman-Wintner law of the iterated
logarithm for martingales.} Ann. Math. Statist. {\bf 41(6)} (1970),
2158-2160.


\bibitem{Wu} L.M. Wu, {\it Forward-backward martingale decomposition and compactness results for additive functionais
of stationary ergodic Markov processes.} Annals de l'i.H.P., S\'erie
Probab and Stat., {\bf 35(2)} (1999), 121-141.
\end{thebibliography}
\end{document}